\documentclass{amsart}

\usepackage{amssymb,amsmath,amsthm,latexsym,booktabs,todonotes}

\theoremstyle{definition}
\newtheorem{definition}{Definition}

\theoremstyle{plain}

\newtheorem{proposition}[definition]{Proposition}

\begin{document}

\title[Canonical forms of Order-$k$ ($k = 2,3,4$) $3 \times \dots \times 3$ Symmetric Tensors]
{Canonical Forms of Order-$k$ ($k=2,3,4$) Symmetric Tensors of Format $3 \times \dots \times 3$ Over Prime Fields}

\author[Stavrou]{Stavros Stavrou}
\address{Department of Mathematics and Statistics, University of Saskatchewan, Canada}
\email{sgs715@mail.usask.ca}





\maketitle

\begin{abstract}
We consider symmetric tensors of format: $3 \times 3$ over $\mathbb{F}_p$ for $p=2,3,5$; $3 \times 3 \times 3$ over $\mathbb{F}_p$ for $p=2,3$; and $3 \times 3 \times 3 \times 3$ over $\mathbb{F}_p$ for $p = 2, 3$. In each case we compute their equivalence classes under the action of the general linear group GL$_3(\mathbb{F}_p)$. We use computer algebra to determine the set of tensors of each symmetric rank, then we compute the orbit of the group action. We determine the maximum symmetric rank of these tensors and compare it with the maximum rank.
\end{abstract}

\section{Introduction}

A tensor is a multidimensional array of numbers. Formally, an order-$k$ tensor is an element of the tensor product of $k$ vector spaces. A first order tensor is a vector, and a second order tensor is a matrix. The rank of a tensor is the minimal number of terms in its expression as a sum of simple tensors. See Cichocki et al. \cite{Cichocki} for a recent survey of matrix and tensor decomposition algorithms, as well as applications in areas such as data mining, email surveillance, gene expression classification, and signal processing. See also Kroonenberg \cite{Kroonenberg}, and Smilde et al. \cite{Smilde}, and references therein.

A tensor is called symmetric if its elements remain fixed under any permutation of the indices. Then we can define the symmetric rank of a symmetric tensor to be the minimal number of terms in its expression as a sum of simple symmetric tensors. Further, we can associate a homogeneous polynomial in $\mathbb{F}[x_0, \dots, x_n]_d$ to any symmetric tensor $X \in S^d V$, and the determination of the symmetric rank of $X$ is equivalent to the Big Waring Problem: determining the minimum integer $r$ such that a generic form of degree $d$ in $n+1$ variables can be written as a sum of $r$ $d$th powers of linear forms \cite{Comon1}. See \cite{B1, B2, Bernardi, B3, B4} and references therein for more information on symmetric tensors, symmetric rank, and decompositions of homogeneous polynomials.

The determination of canonical forms of tensors is a generalization of the matrix case. For any $m \times n$ matrix $A$ (with $m \leq n$) there exist invertible matrices $P \in GL_m(\mathbb{F})$ and $Q \in GL_n(\mathbb{F})$ such that 
\[
PAQ = 
\begin{bmatrix} I_r & 0 \\ 0 & 0 \end{bmatrix},
\]
where $I_r$ is the $r \times r$ identity and $r \leq m$ is the rank. The matrix $PAQ$ is called the Smith normal form (canonical form), and is the orbit representative for the action of $GL_m(\mathbb{F}) \times GL_n(\mathbb{F})$, the direct product of general linear groups. When $m, n = 2$, there are three canonical forms,
\[
\begin{bmatrix} 0 & 0 \\ 0 & 0 \end{bmatrix}, \qquad
\begin{bmatrix} 1 & 0 \\ 0 & 0 \end{bmatrix}, \qquad
\begin{bmatrix} 1 & 0 \\ 0 & 1 \end{bmatrix}.
\]
Determining the canonical forms of order-$k$ tensors for $k>2$ is a more complicated problem. For the next simplest case, consider $2 \times 2 \times 2$ tensors. Le Paige (1881) determined there are 7 canonical forms over $\mathbb{C}$ \cite{LePaige}, and later Oldenburger determined there are 8 over $\mathbb{R}$ \cite{Oldenburger1}. In previous work, we considered the same problem over the prime fields $\mathbb{F}_2$ and $\mathbb{F}_3$ \cite{BremnerStavrou}. The maximum rank of these tensors over $\mathbb{R}$, $\mathbb{C}$, and the prime fields we considered is 3. For larger tensor formats over $\mathbb{R}$ and $\mathbb{C}$, the problem is more complicated and, in general, there is not a finite classification \cite{deSilvaLim}. 

In previous work, we considered the same problem of determining canonical forms over prime fields, but restricted our attention to symmetric tensors. For the $2 \times 2 \times 2$ and $2 \times 2 \times 2 \times 2$ symmetric formats over $\mathbb{R}$ and $\mathbb{C}$, this problem was analyzed independently by Gurevich \cite{Gurevich} first, then by Weinberg \cite{Weinberg}. 
In \cite{Stavrou}, we extended these results by determining the canonical forms of $2 \times 2 \times 2$ symmetric tensors over the prime fields $\mathbb{F}_p$ for $p = 2, 3, 5, 7, 11, 13, 17$, as well as $2 \times 2 \times 2 \times 2$ symmetric tensors over $\mathbb{F}_p$ for $p = 2, 3, 5, 7$. 

We continue in this direction by considering order-$k$ (for $k = 2,3,4$) symmetric tensors of format $3 \times \dots \times 3$ over $\mathbb{F}_p$, for $p \in \{2,3,5\}$. For $3 \times 3$ symmetric tensors (i.e. symmetric matrices), we work over $\mathbb{F}_p$ for $p = 2, 3, 5$. The maximum symmetric rank is 4, which is larger than the maximum rank of 3.  For $3 \times 3 \times 3$ symmetric tensors, we work over $\mathbb{F}_p$ for $p=2,3$. The maximum symmetric rank is 7, which is larger than the maximum rank of 5 over $\mathbb{C}$ \cite{PierreComon}. Lastly, for $3 \times 3 \times 3 \times 3$ symmetric tensors, we work over $\mathbb{F}_p$ for $p=2, 3$. Over $\mathbb{F}_3$, the number of symmetric tensors was too large to determine, and so the procedure had to be terminated. At termination, there were symmetric tensors with symmetric rank 13, thus we have a lower bound on the symmetric rank.

In each case, we compute the symmetric tensors of each rank. The general linear group GL$_3(\mathbb{F}_p)$ partitions the tensors of each rank into a disjoint union of orbits, where the tensors in each orbit are equivalent. The elements are arranged in lexical order and the minimal element is the canonical form of its orbit under this group action. We comment here that due to limited computer memory, we are not able to consider larger prime fields, or other tensor formats.

\section{Preliminaries}

An order-$k$ {\bf tensor} $X$ is an element of the {\bf tensor product} of $k$ vector spaces $V_1 \otimes \dots \otimes V_k$, where the {\bf order} refers to the number of dimensions. Once we fix a basis in each vector space, we can associate to $X$ a $k$-dimensional array. Then an order-1 tensor is a vector, an order-2 tensor is a matrix, and we denote an order-$k$ tensor of format $d_1 \times \dots \times d_k$ as $[x_{i_1 \dots i_k}] \in \mathbb{F}^{d_1 \times \dots \times d_k}$. We can represent the tensor $X$ in {\bf vectorized (flattened)} form by writing the columns of $X$ in a vector format, where the entries are in lexical order of the $k$-tuples of the subscripts. 

An order-$k$ tensor $[x_{i_1 \dots i_k}] \in \mathbb{F}^{n \times \dots \times n}$ is called {\bf symmetric} if 
\[
x_{i_{\pi(1)} \dots i_{\pi(k)}} = x_{i_1 \dots i_k}, \qquad i_1, \dots, i_k \in \{1, \dots, n\},
\]
for all permutations $\pi \in S_k$  \cite{Comon1}. 

\begin{proposition} [Proposition 3.7 \cite{Comon1}]
Let $X = [x_{i_1 \dots i_k}] \in \mathbb{F}^{n \times \dots \times n}$ be an order-$k$ tensor. Then 
\[
\pi(X) = X
\]
for all permutations $\pi \in S_k$ if and only if 
\[
x_{i_{\pi(1)} \dots i_{\pi(k)}} = x_{i_1 \dots i_k}, \qquad i_1, \dots, i_k \in \{1, \dots, n\}
\]
for all permutations $\pi \in S_k$.
\end{proposition}

In this paper we consider $3 \times \dots \times 3$ symmetric tensors of order-$k$ for $k=2,3,4$. For $k = 2$, a symmetric matrix has the form
\[
\begin{bmatrix}
a & b & c \\
b & d & e \\
c & e & f
\end{bmatrix},
\]
where the entries are scalars in $\mathbb{F}_p$. For $k = 3$, 
a $3 \times 3 \times 3$ symmetric tensor, in terms of its $3 \times 3$ frontal slices, has the form 
\[
\left[
\begin{array}{ccc|ccc|ccc}
a & b & c & b & d & e & c & e & f \\
b & d & e & d & j & g & e & g & h \\
c & e & f & e & g & h & f & h & k
\end{array}
\right].
\]
Lastly, for $k = 4$, a $3 \times 3 \times 3 \times 3$ symmetric tensor, in terms of its $3 \times 3$ frontal slices, has the form
\[
\left[
\begin{array}{ccc|ccc|ccc}
a & b & f & b & c & g & f & g & j \\
b & c & g & c & d & i & g & i & h \\
f & g & j & g & i & h & j & h & k \\ \hline
b & c & g & c & d & i & g & i & h \\
c & d & i & d & e & l & i & l & m \\
g & i & h & i & l & m & h & m & n \\ \hline
f & g & j & g & i & h & j & h & k \\
g & i & h & i & l & m & h & m & n \\
j & h & k & h & m & n & k & n & p
\end{array}
\right].
\]

We will denote the set of all order-$k$, $n$-dimensional symmetric tensors over the field $\mathbb{F}$ by $\mathcal{S}^k(\mathbb{F}^n) \subset \mathbb{F}^{n \times \dots \times n}$. The set of such tensors satisfies the property $\pi(X) = X$ for all $\pi \in S_k$ and $X \in \mathcal{S}^k(\mathbb{F}^n)$. 

\begin{definition}\cite{Hitchcock1} 
A tensor $X \in \mathbb{F}^{d_1 \times \dots \times d_k}$ is {\bf simple} if it can be written as
\[
X = u^{(1)} \otimes u^{(2)} \otimes \cdots \otimes u^{(k)},
\]
with non-zero $u^{(i)} \in \mathbb{F}^{d_i}$ for $i = 1, \dots, k$. The $(i_1, i_2, \dots, i_k)$th entry of $X$ is 
\[
x_{ i_1 i_2 \cdots i_k } = u^{(1)}_{i_1} u^{(2)}_{i_2} \cdots u^{(k)}_{i_k}.
\]
\end{definition}

\begin{definition}
A tensor $X \in \mathcal{S}^k(\mathbb{F}^n)$ is called {\bf simple symmetric} if it can be written as 
\[
X = u \otimes \dots \otimes u
\]
with $k$-many non-zero vectors $u \in \mathbb{F}^n$. The $(i_1, \dots, i_k)$th entry of $X$ is 
\[
x_{i_1 \dots i_k} = u_{i_1} \dots u_{i_k}.
\]
\end{definition}

\begin{definition}\cite{Hitchcock1}
A tensor has {\bf outer product rank} $r$ if it can be written as a sum of $r$ (and no fewer) decomposable tensors, 
\[
X = \sum_{i = 1}^r  u_i^{(1)} \otimes \dots \otimes u_i^{(k)} = u_1^{(1)} \otimes \dots \otimes u_1^{(k)} + \dots + u_r^{(1)} \otimes \dots \otimes u_r^{(k)}
\]
where $u_i^{(1)} \in \mathbb{F}^{d_1}, \dots, u_i^{(k)} \in \mathbb{F}^{d_k}, i = 1, \dots, r$.
We write rank$(X)$ to denote the outer product rank of $X$.
\end{definition}

\begin{definition}
A tensor has {\bf symmetric outer product rank} $s$ if it can be written as a sum of $s$ (and no fewer) symmetric simple tensors,
\[
X = \sum_{i=1}^s u_i^{\otimes k}.
\]
We write rank$_S(X)$ to denote the symmetric outer product rank of $X$. 
\end{definition}
The only rank-0 symmetric tensor is the zero tensor. We will drop the words outer product and simply say symmetric rank. 

\begin{definition}
The {\bf maximum rank} is defined to be  
\[
\text{max}\{ \text{rank}(X) \mid X \in \mathbb{F}^{d_1 \times \dots \times d_k} \}. 
\]
If we replace rank$(X)$ with rank$_S(X)$ we get the analogous definition of maximum symmetric rank.

\end{definition}

In order to compute the canonical forms of symmetric tensors we use a tensor-matrix multiplication. {\bf Multilinear matrix multiplication} is a tensor-matrix multiplication that allows us to multiply matrices on each of the modes of the tensor. If $X = [x_{i_1 \dots i_k}] \in \mathbb{F}^{d_1 \times \dots \times d_k}$ and 
\[
A_1 = [a^{(1)}_{u_1 i_1}] \in \mathbb{F}^{c_1 \times d_1}, \dots, A_k = [a^{(k)}_{u_k i_k}] \in \mathbb{F}^{c_k \times d_k},
\]
then $Y= (A_1, \dots, A_k)\cdot X = [y_{u_1 \dots u_k}] \in \mathbb{F}^{c_1 \times \dots \times c_k}$ is the new tensor defined by
\[
y_{u_1 \dots u_k} = \sum_{i_1, \dots, i_k = 1}^{d_1, \dots, d_k} a^{(1)}_{u_1 i_1}  \dots  a^{(k)}_{u_k i_k} x_{i_1 \dots i_k}.
\]
Since we are restricting ourselves to symmetric tensors, we impose the condition that $A_1 = \dots = A_k$, otherwise, multilinear matrix multiplication can transform a symmetric tensor into a non-symmetric tensor.


We consider the action of the symmetry group GL$_3(\mathbb{F}_p)$, which does not change the rank of a tensor \cite{deSilvaLim}. The action of this group decomposes the set of order-$k$ symmetric tensors into a disjoint union of orbits, where the tensors in each orbit are equivalent under the group action. The orbit of $X$ is the set 
\[
\mathcal{O}_X := \{ (g, \dots, g) \cdot X \mid g \in \text{GL}_3(\mathbb{F}_p)\}.
\] 
We define the {\bf canonical form} of $X$ to be the minimal element in its orbit with respect to the lexical ordering.

\section{Algorithms}

Using computer algebra, we generate all the symmetric tensors in each symmetric rank. Then we apply the group action which decomposes the set of tensors within each rank into a disjoint union of orbits. The minimal element of each orbit with respect to the lexical ordering is the canonical form.  

We determine the symmetric tensors in each symmetric rank by first generating the set $R_1$ of simple symmetric tensors, which is achieved by computing all possible products $v \otimes \dots \otimes v$ of non-zero vectors $v \in \mathbb{F}^3_p$, for each prime. Then to determine the elements in $R_i$ we compute all possible sums $(T + S)\ \textrm{mod} \ p$ for $T \in R_{i-1}$ and $S \in R_1$, and then subtract the tensors already computed in the previous ranks: $R_i = R_i \setminus \cup_{k=1}^{i-1} R_k$. The procedure terminates once no new symmetric tensors are generated when taking sums.  In our previous work \cite{Stavrou}, when $T+S$ was computed, the procedure searched to see if $T+S$ already existed in the current rank and the lower ranks. If it did not, then it was added to the current rank set. In this paper, we made the necessary modification of subtracting the sets so that the ranks were computed more efficiently. This is what allowed us to consider larger formats than what we considered in \cite{Stavrou}.


Next we apply the action of the general linear group to decompose the tensors in each rank into a disjoint union of orbits. The pseudocode is provided in Table \ref{groupalgorithmtable}. This algorithm was used in our previous paper \cite{Stavrou}. In this paper, we made two major modifications to the program. The first modification was removing an index (mode) from the procedures in order to consider second order $3 \times 3$ symmetric tensors.  The second modification was to index the elements to 3, since each mode has length 3. The pseudocode in Table \ref{groupalgorithmtable} displays the order-4 case. For the order-3 case, remove the $\ell$ index, and for the order-2 case, remove the $\ell$ and $k$ indices.

\begin{table}
\begin{itemize}
\item[]
\texttt{unflatten}$( x )$
\item[] \quad
set $t \leftarrow 0$
\item[] \quad
for $i, j, k, \ell = 1, 2, 3$ do:
set $t \leftarrow t + 1$;
set $y_{ijk\ell} \leftarrow x_t$
\item[] \quad
$\texttt{return}( y )$
\end{itemize}
\medskip
\begin{itemize}
\item[]
\texttt{groupaction}$( g, x, m )$
\item[] \quad
set $y \leftarrow \texttt{unflatten}( x )$
\item[] \quad
if $m = 1$ then
for $j,k,\ell = 1,2,3$ do:
\item[] \quad \quad
set $v \leftarrow [ \, y_{1jk\ell}, \, y_{2jk\ell} \, ]$;
\item[] \quad \quad
set $w \leftarrow [ \, g_{11} v_1 {+} g_{12} v_2 {+} g_{13}v_3 \, \mathrm{mod}\,p, 
\dots,
\, g_{31} v_1 {+} g_{32} v_2{+} g_{33} v_3 \, \mathrm{mod}\,p \,]$
\item[] \quad \quad
for $i = 1, 2, 3$ do: set $y_{ijk\ell} \leftarrow w_i$
\item[] \quad
if $m = 2$ then \dots \emph{(similar for second subscript)}
\item[] \quad
if $m = 3$ then \dots \emph{(similar for third subscript)}
\item[] \quad
if $m = 4$ then \dots \emph{(similar for fourth subscript)}
\item[] \quad
\texttt{return}( \texttt{flatten}( $y$ ) )
\end{itemize}
\medskip
\begin{itemize}
\item[]
\texttt{smallorbit}$( x )$
\item[] \quad
set $\texttt{result} \leftarrow \{\,\}$
\item[] \quad
for $a \in GL_3(\mathbb{F}_p)$ do:
\item[] \quad \quad
set $y \leftarrow \texttt{groupaction}( a, x, 1 )$
\item[] \quad \quad
set $z \leftarrow \texttt{groupaction}( a, y, 2 )$
\item[] \quad \quad 
set $w \leftarrow \texttt{groupaction}( a, z, 3 )$
\item[] \quad \quad 
set $u \leftarrow \texttt{groupaction}( a, w, 4 )$
\item[] \quad
set $\texttt{result} \leftarrow \texttt{result} \cup \{ u \}$
\item[] \quad
\texttt{return}( \texttt{result} )
\end{itemize}
\medskip
\begin{itemize}
\item
for $r = 0,\dots,\texttt{maximumrank}$ do:
\item[] \quad
set $\texttt{representatives}[r] \leftarrow \{\,\}$;
set $\texttt{remaining} \leftarrow \texttt{arrayset}[r]$
\item[] \quad
while $\texttt{remaining} \ne \{\,\}$ do:
\item[] \quad \quad
set $x \leftarrow \texttt{remaining}[1]$;
set $\texttt{xorbit} \leftarrow \texttt{largeorbit}( x )$
\item[] \quad \quad
append $\texttt{xorbit}[1]$ to $\texttt{representatives}[r]$
\item[] \quad \quad
set $\texttt{remaining} \leftarrow \texttt{remaining} \setminus \texttt{xorbit}$
\end{itemize}
\medskip
\caption{Algorithm for group action (pseudocode)}
\label{groupalgorithmtable}
\end{table}

\section{Symmetric Tensors of Format $3 \times 3$}

Every order-$k$ symmetric tensor of dimension $n$ may be uniquely associated with a homogeneous polynomial (also called a quantic) of degree $k$ in $n$ variables \cite{Comon1}. The problem of determining symmetric rank is equivalent to the Big Waring Problem: determining the minimal number of $p$th powers of linear terms \cite{Ehrenborg2} \cite{Ehrenborg3}. The Alexander-Hirschowitz (AH) Theorem \cite{AH} completely describes the calculation of the generic rank of symmetric tensors.

\subsection{Canonical Forms of $3 \times 3$ Symmetric Tensors over $F_2$.}

There are 512 $3 \times 3$ tensors over $\mathbb{F}_2$, where 64 are symmetric. Every symmetric tensor has a symmetric decomposition. The maximum symmetric rank is 3, which is equal to the maximum rank. The symmetric ranks, orders of each orbit, and the minimal representatives of each orbit are given in Table \ref{table33symmetricmod2}. The number of tensors in each symmetric rank is listed below.
\[
\begin{array}{lrrrrrrrr}
\text{rank} & 0 & 1 & 2 & 3 \\
\text{number} & 1 & 7 & 21 & 35 \\
\text{$\approx$ $\%$} & 1.5625\% & 10.9375\% & 32.8125\% & 54.6875\%
\end{array}
\]

  \begin{table}
  \[
  \begin{array}{ccc}
  \text{symmetric rank} & \text{orbit size} &  \text{canonical form} 
    \\
  \toprule 
0 & 1 & \begin{bmatrix}0&0&0\\0&0&0\\0&0&0\end{bmatrix}
\\
\midrule 
1 & 7 & \begin{bmatrix}0&0&0\\0&0&0\\0&0&1\end{bmatrix}
\\
\midrule
2 & 21 &  \begin{bmatrix}0&0&0\\0&0&1\\0&1&1\end{bmatrix}
\\
\midrule
3 & 7 & \begin{bmatrix}0&0&0\\0&0&1\\0&1&0\end{bmatrix}
\\
\midrule
3 & 28 & \begin{bmatrix}0&0&1\\0&1&0\\1&0&0\end{bmatrix}
  \\
  \bottomrule
  \end{array}
  \]
  \medskip
  \caption{Canonical forms of $3 \times 3$ symmetric tensors over $\mathbb{F}_2$}
  \label{table33symmetricmod2}
  \end{table}

\subsection{Canonical Forms of $3 \times 3$ Symmetric Tensors over $F_3$.}

There are 19,683 $3 \times 3$ tensors over $\mathbb{F}_3$, where 729 are symmetric. Every symmetric tensor has a symmetric decomposition. The maximum symmetric rank in this case is 4, which is larger than the maximum rank of 3. The symmetric ranks, orders of each orbit, and the minimal representatives of each orbit are given in Table \ref{table33symmetricmod3}. The number of tensors in each symmetric rank is listed below.
\[
\begin{array}{lrrrrrrrr}
\text{rank} & 0 & 1 & 2 & 3 & 4 \\
\text{number} & 1 & 13 & 91 & 390 & 234 \\
\text{$\approx$ $\%$} & 0.1372\% & 1.7833\% & 12.4829\% & 53.4979\% & 32.0988
\end{array}
\]

  \begin{table}
  \[
  \begin{array}{ccc}
  \text{symmetric rank} & \text{orbit size} &  \text{canonical form} 
    \\
  \toprule 
0 & 1 & \begin{bmatrix}0&0&0\\0&0&0\\0&0&0\end{bmatrix}
\\
\midrule 
1 & 13 & \begin{bmatrix}0&0&0\\0&0&0\\0&0&1\end{bmatrix}
\\
\midrule
2 & 13 &  \begin{bmatrix}0&0&0\\0&0&0\\0&0&2\end{bmatrix}
\\
\midrule
2 & 78 &  \begin{bmatrix}0&0&0\\0&1&0\\0&0&1\end{bmatrix}
\\
\midrule
3 & 156 & \begin{bmatrix}0&0&0\\0&0&1\\0&1&0\end{bmatrix}
\\
\midrule
3 & 234 &  \begin{bmatrix}0&0&1\\0&2&0\\1&0&0\end{bmatrix}
\\
\midrule
4 & 234 & \begin{bmatrix}0&0&1\\0&1&0\\1&0&0\end{bmatrix}
  \\
  \bottomrule
  \end{array}
  \]
  \medskip
  \caption{Canonical forms of $3 \times 3$ symmetric tensors over $\mathbb{F}_3$}
  \label{table33symmetricmod3}
  \end{table}

\subsection{Canonical Forms of $3 \times 3$ Symmetric Tensors over $F_5$.}

There are 1,953,125 $3 \times 3$ tensors over $\mathbb{F}_5$, where 15,625 are symmetric.  The maximum symmetric rank is again 4. The order of the group GL$(\mathbb{F}_5)$ is 427,307 which our program cannot handle, and thus we do not determine the orbits under the group action. The symmetric ranks, number of tensors in each rank, and the minimal representative of each rank are given in Table \ref{table33symmetricmod5}. 
\[
\begin{array}{lrrrrrrrr}
\text{rank} & 0 & 1 & 2 & 3 & 4 \\
\text{number} & 1 & 62 & 1922 & 7440 & 6200 \\
\text{$\approx$ $\%$} & 0.0064\% & 0.3968\% & 12.3008\% & 47.6160\% & 39.6800
\end{array}
\]

  \begin{table}
  \[
  \begin{array}{ccc}
  \text{symmetric rank} & \text{rank size} &  \text{minimal element} 
    \\
  \toprule 
0 & 1 & \begin{bmatrix}0&0&0\\0&0&0\\0&0&0\end{bmatrix}
\\
\midrule 
1 & 62 & \begin{bmatrix}0&0&0\\0&0&0\\0&0&1\end{bmatrix}
\\
\midrule
2 & 1922 &  \begin{bmatrix}0&0&0\\0&0&0\\0&0&2\end{bmatrix}
\\
\midrule
3 & 7440 &  \begin{bmatrix}0&0&0\\0&1&0\\0&0&2\end{bmatrix}
\\
\midrule
4 & 6200 & \begin{bmatrix}0&0&1\\0&2&0\\1&0&0\end{bmatrix}
  \\
  \bottomrule
  \end{array}
  \]
  \medskip
  \caption{Minimal elements of $3 \times 3$ symmetric tensors over $\mathbb{F}_5$}
  \label{table33symmetricmod5}
  \end{table}

\section{Symmetric Tensors of Format $3 \times 3 \times 3$}

In Dickson's 1908 paper \cite{Dickson}, he cites that a complete set of canonical forms of ternary cubics was determined over $\mathbb{C}$. He considers canonical forms by determining the algebraic irrationalities occurring in the reducing linear transformations over $\mathbb{F}_3$, and finds 11 distinct forms.

More recently, Comon \cite{PierreComon} summarize the generic ranks of symmetric tensors of some specific formats. In particular, they tabulate the equivalence classes for ternary cubics (i.e. symmetric $3 \times 3 \times 3$ tensors). Below are the orbits under the action of the group of invertible three-dimensional changes of coordinates.
\[
\begin{array}{c|c}
\text{Orbit} & \text{rank} \\ \hline
x^3 & 1 \\
x^3+y^3 & 2 \\
x^2y & 3 \\
x^3 + 3y^2z & 4 \\
x^3 + y^3 + 6xyz & 4 \\
x^3 + 6xyz & 4 \\
a(x^3 + y^3 + z^3) + 6bxyz & 4 \\
x^2y + xz^2 & 5
\end{array}
\]

Kogan and Maza \cite{Kogan} determined the equivalence classes of ternary cubics under general complex linear changes of variables using a computational approach. See in particular Theorem 5 of \cite{Kogan}. Their contribution provided a method of computing the signature manifolds for each of the equivalence classes, and their algorithm matches a cubic with its canonical form while producing the required linear transformation explicitly.

\subsection{Canonical Forms of $3 \times 3 \times 3$ Symmetric Tensors over $\mathbb{F}_2$}

There are 1024 $3 \times 3 \times 3$ symmetric tensors over $\mathbb{F}_2$. Only 128 (12.5\%) of these have a symmetric decomposition. The maximum symmetric rank is 7. The symmetric ranks, orders of each orbit, and the minimal representatives (in flattened form) of each orbit are given in Table \ref{table333symmetricmod2}. To make the canonical forms more legible we replace the zero entries with $\cdot$. The number of tensors in each symmetric rank is listed below (percentages add to 12.5\%).
\[
\begin{array}{lrrrrrrrr}
\text{rank} & 0 & 1 & 2 & 3 & 4 & 5 & 6 & 7 \\
\text{number} & 1 & 7 & 21 & 35 & 35 & 21 & 7 & 1 \\
\text{$\approx$ $\%$} & 0.10\% & 0.68\% & 2.05\% & 3.42\% & 3.42\% & 2.05\% & 0.68\% & 0.10 \%
\end{array}
\]

  \begin{table}
  \[
  \begin{array}{ccc}
  \text{symmetric rank} & \text{orbit size} &  \text{canonical form (flattened)} 
    \\
  \toprule 
0 & 1 & [ \cdot  \cdot  \cdot  \cdot  \cdot  \cdot  \cdot  \cdot  \cdot  \cdot  \cdot  \cdot  \cdot  \cdot  \cdot  \cdot  \cdot  \cdot  \cdot  \cdot  \cdot  \cdot  \cdot  \cdot  \cdot  \cdot  \cdot ]
\\
\midrule 
1 & 7 & [\cdot \cdot \cdot \cdot \cdot \cdot \cdot \cdot \cdot \cdot \cdot \cdot \cdot \cdot \cdot \cdot \cdot \cdot \cdot \cdot \cdot \cdot \cdot \cdot \cdot \cdot 1] 
\\
\midrule
2 & 21 &  [\cdot \cdot \cdot \cdot \cdot \cdot \cdot \cdot \cdot \cdot \cdot \cdot \cdot \cdot 1 \cdot 1 1 \cdot \cdot \cdot \cdot 1 1 \cdot 1 1]
\\
\midrule
3 & 7 & [\cdot \cdot \cdot \cdot \cdot \cdot \cdot \cdot \cdot \cdot \cdot \cdot \cdot \cdot 1 \cdot 1 1 \cdot \cdot \cdot \cdot 1 1 \cdot 1 \cdot]
\\
3 & 18 &  [\cdot \cdot 1 \cdot \cdot \cdot 1 \cdot 1 \cdot \cdot \cdot \cdot 1 \cdot \cdot \cdot \cdot 1 \cdot 1 \cdot \cdot \cdot 1 \cdot 1]
\\
\midrule
4 & 7 &  [\cdot \cdot \cdot \cdot \cdot 1 \cdot 1 \cdot \cdot \cdot 1 \cdot \cdot 1 1 1 1 \cdot 1 \cdot 1 1 1 \cdot 1 \cdot]  
\\  
4 & 18 & [\cdot \cdot 1 \cdot \cdot \cdot 1 \cdot 1 \cdot \cdot \cdot \cdot 1 \cdot \cdot \cdot \cdot 1 \cdot 1 \cdot \cdot \cdot 1 \cdot \cdot]
\\   
\midrule
5 & 21 & [\cdot \cdot \cdot \cdot \cdot 1 \cdot 1 \cdot \cdot \cdot 1 \cdot \cdot 1 1 1 1 \cdot 1 \cdot 1 1 1 \cdot 1 1]    
\\
\midrule 
6 & 7 & [\cdot \cdot \cdot \cdot \cdot 1 \cdot 1 \cdot \cdot \cdot 1 \cdot \cdot \cdot 1 \cdot \cdot \cdot 1 \cdot 1 \cdot \cdot \cdot \cdot 1]
\\  
\midrule
7 & 1 & [\cdot \cdot \cdot \cdot \cdot 1 \cdot 1 \cdot \cdot \cdot 1 \cdot \cdot \cdot 1 \cdot \cdot \cdot 1 \cdot 1 \cdot \cdot \cdot \cdot \cdot]   
  \\
  \bottomrule
  \end{array}
  \]
  \medskip
  \caption{Canonical forms of $3 \times 3 \times 3$ symmetric tensors over $\mathbb{F}_2$}
  \label{table333symmetricmod2}
  \end{table}

\subsection{Canonical Forms of $3 \times 3 \times 3$ Symmetric Tensors over $\mathbb{F}_3$}

There are 59,049 $3 \times 3 \times 3$ symmetric tensors over $\mathbb{F}_3$, and unlike in the previous case, every symmetric tensor has a symmetric decomposition. The maximum symmetric rank is 7. The symmetric ranks, orders of each orbit, and the minimal representatives of each orbit are given in Table \ref{table333symmetricmod3}. The number of tensors in each rank is listed in the table below. We mention now that we cannot consider larger primes because of insufficient computer memory. There are over 9.7 million symmetric $3 \times 3 \times 3$ tensors over $\mathbb{F}_5$.

\[
\begin{array}{lrrrrrrrr}
\text{rank} & 0 & 1 & 2 & 3 & 4 & 5 & 6 & 7 \\
\text{number} & 1 & 26 & 312 & 2288 & 11440 & 30342 & 14352 & 288 \\
\text{$\approx$ $\%$} & 0.00\% & 0.05\% & 0.53\% & 3.87\% & 19.37\% & 51.38\% & 24.31\% & 0.49 \%
\end{array}
\]

\begin{table}
 \[
  \begin{array}{cccc}
  \text{symmetric rank} & \text{orbit size} &  \text{canonical form (flattened)} 
    \\
  \toprule
0 & 1 &  [\cdot  \cdot  \cdot  \cdot  \cdot  \cdot  \cdot  \cdot  \cdot  \cdot  \cdot  \cdot  \cdot  \cdot  \cdot  \cdot  \cdot  \cdot  \cdot  \cdot  \cdot  \cdot  \cdot  \cdot  \cdot  \cdot  \cdot]   
\\
\midrule  
1 & 26 & [\cdot \cdot \cdot \cdot \cdot \cdot \cdot \cdot \cdot \cdot \cdot \cdot \cdot \cdot \cdot \cdot \cdot \cdot \cdot \cdot \cdot \cdot \cdot \cdot \cdot \cdot 1]  
\\
\midrule
2 & 312 &  [\cdot \cdot \cdot \cdot \cdot \cdot \cdot \cdot \cdot \cdot \cdot \cdot \cdot \cdot 1 \cdot 1 \cdot \cdot \cdot \cdot \cdot 1 \cdot \cdot \cdot 1]  
\\
\midrule  
3 & 104 &  [\cdot \cdot \cdot \cdot \cdot \cdot \cdot \cdot \cdot \cdot \cdot \cdot \cdot \cdot \cdot \cdot \cdot 1 \cdot \cdot \cdot \cdot \cdot 1 \cdot 1 \cdot]  
\\ 
3 & 312 &  [\cdot \cdot \cdot \cdot \cdot \cdot \cdot \cdot \cdot \cdot \cdot \cdot \cdot \cdot 1 \cdot 1 \cdot \cdot \cdot \cdot \cdot 1 \cdot \cdot \cdot 2]    
\\
3 & 1872 & [\cdot \cdot 1 \cdot \cdot \cdot 1 \cdot \cdot \cdot \cdot \cdot \cdot 1 \cdot \cdot \cdot \cdot 1 \cdot \cdot \cdot \cdot \cdot \cdot \cdot 1]     
\\
\midrule  
4 & 208 & [\cdot \cdot \cdot \cdot \cdot \cdot \cdot \cdot \cdot \cdot \cdot \cdot \cdot \cdot \cdot \cdot \cdot 1 \cdot \cdot \cdot \cdot \cdot 1 \cdot 1 1]  
\\
4 & 1872 & [\cdot \cdot \cdot \cdot \cdot \cdot \cdot \cdot 1 \cdot \cdot \cdot \cdot 1 \cdot \cdot \cdot \cdot \cdot \cdot 1 \cdot \cdot \cdot 1 \cdot \cdot ]    
\\
4 & 468 & [\cdot \cdot \cdot \cdot \cdot 1 \cdot 1 \cdot \cdot \cdot 1 \cdot \cdot \cdot 1 \cdot \cdot \cdot 1 \cdot 1 \cdot \cdot \cdot \cdot \cdot]    
\\
4 & 1404 & [\cdot \cdot 1 \cdot \cdot \cdot 1 \cdot \cdot \cdot \cdot \cdot \cdot \cdot 1 \cdot 1 \cdot 1 \cdot \cdot \cdot 1 \cdot \cdot \cdot 1]    
\\
4 & 5616 & [\cdot \cdot 1 \cdot \cdot \cdot 1 \cdot \cdot \cdot \cdot \cdot \cdot 1 \cdot \cdot \cdot \cdot 1 \cdot \cdot \cdot \cdot \cdot \cdot \cdot 2]    
\\
4 & 1872 & [\cdot \cdot 1 \cdot \cdot \cdot 1 \cdot \cdot \cdot \cdot \cdot \cdot 1 \cdot \cdot \cdot 1 1 \cdot \cdot \cdot \cdot 1 \cdot 1 1]   
\\
\midrule  
5 & 624 & [\cdot \cdot \cdot \cdot \cdot \cdot \cdot \cdot 1 \cdot \cdot \cdot \cdot \cdot 1 \cdot 1 \cdot \cdot \cdot 1 \cdot 1 \cdot 1 \cdot 1]     
\\
5 & 3744 & [\cdot \cdot \cdot \cdot \cdot \cdot \cdot \cdot 1 \cdot \cdot \cdot \cdot 1 \cdot \cdot \cdot \cdot \cdot \cdot 1 \cdot \cdot \cdot 1 \cdot 1]    
\\
5 & 2808  & [\cdot \cdot \cdot \cdot \cdot 1 \cdot 1 \cdot \cdot \cdot 1 \cdot \cdot \cdot 1 \cdot \cdot \cdot 1 \cdot 1 \cdot \cdot \cdot \cdot 1]    
\\
5 & 5616  &   [\cdot \cdot \cdot \cdot \cdot 1 \cdot 1 \cdot \cdot \cdot 1 \cdot 1 \cdot 1 \cdot \cdot \cdot 1 \cdot 1 \cdot \cdot \cdot \cdot 1]    
\\
5 & 702  & [\cdot \cdot \cdot \cdot 1 \cdot \cdot \cdot 1 \cdot 1 \cdot 1 \cdot \cdot \cdot \cdot \cdot \cdot \cdot 1 \cdot \cdot \cdot 1 \cdot \cdot]    
 \\
5 & 5616 & [\cdot \cdot 1 \cdot \cdot \cdot 1 \cdot \cdot \cdot \cdot \cdot \cdot 1 \cdot \cdot \cdot 1 1 \cdot \cdot \cdot \cdot 1 \cdot 1 \cdot]    
\\ 
5 & 5616  &  [\cdot \cdot 1 \cdot \cdot \cdot 1 \cdot \cdot \cdot \cdot \cdot \cdot 1 \cdot \cdot \cdot 1 1 \cdot \cdot \cdot \cdot 1 \cdot 1 2]    
\\
5 & 5616  & [\cdot \cdot 1 \cdot \cdot \cdot 1 \cdot \cdot \cdot \cdot \cdot \cdot 1 \cdot \cdot \cdot 2 1 \cdot \cdot \cdot \cdot 2 \cdot 2 \cdot]    
\\
\midrule  
6 & 624 & [\cdot \cdot \cdot \cdot \cdot \cdot \cdot \cdot 1 \cdot \cdot \cdot \cdot \cdot 1 \cdot 1 \cdot \cdot \cdot 1 \cdot 1 \cdot 1 \cdot \cdot]   
\\
6 & 624 & [\cdot \cdot \cdot \cdot \cdot \cdot \cdot \cdot 1 \cdot \cdot \cdot \cdot \cdot 1 \cdot 1 \cdot \cdot \cdot 1 \cdot 1 \cdot 1 \cdot 2]     
\\ 
6 & 5616 & [\cdot \cdot \cdot \cdot 1 \cdot \cdot \cdot 1 \cdot 1 \cdot 1 \cdot \cdot \cdot \cdot \cdot \cdot \cdot 1 \cdot \cdot \cdot 1 \cdot 1]    
\\
6 & 3744 & [\cdot \cdot 1 \cdot 1 \cdot 1 \cdot \cdot \cdot 1 \cdot 1 \cdot 2 \cdot 2 1 1 \cdot \cdot \cdot 2 1 \cdot 1 1]    
\\
6 & 3744 & [\cdot \cdot 1 \cdot 1 \cdot 1 \cdot \cdot \cdot 1 \cdot 1 \cdot 2 \cdot 2 1 1 \cdot \cdot \cdot 2 1 \cdot 1 2]    
 \\
 \midrule 
7 & 288  & [\cdot \cdot 1 \cdot 1 \cdot 1 \cdot \cdot \cdot 1 \cdot 1 \cdot 2 \cdot 2 1 1 \cdot \cdot \cdot 2 1 \cdot 1 \cdot] 
\\ 
  \bottomrule
 \end{array}
  \]
  \medskip
  \caption{Canonical forms of $3 \times 3 \times 3$ symmetric tensors over $\mathbb{F}_3$}
  \label{table333symmetricmod3}
  \end{table}

\section{Symmetric Tensors of Format $3 \times 3 \times 3 \times 3$}

We consider $3 \times 3 \times 3 \times 3$ symmetric tensors over $\mathbb{F}_2$. Since the order of the group is too large for our algorithm, we do not determine the orbits under the group action. Instead, we determine the maximum symmetric rank, the number of symmetric tensors in each rank, and the minimal element.

\subsection{Canonical Forms of $3 \times 3 \times 3 \times 3$ Symmetric Tensors over $\mathbb{F}_2$}
There are $2^{81}$ tensors over $\mathbb{F}_2$, where 32,768 are symmetric. Only 128 (approximately 0.29\%) of these symmetric tensors have a symmetric decomposition.
\[
\begin{array}{lrrrrrrrr}
\text{rank} & 0 & 1 & 2 & 3 & 4 & 5 & 6 & 7 \\
\text{number} & 1 & 7 & 21 & 35 & 35 & 21 & 7 & 1 \\
\text{$\approx$ $\%$} & 0.0031\% & 0.0214\% & 0.0641\% & 0.1068\% & 0.1068\% & 0.0641\% & 0.0214\% & 0.0031 \%
\end{array}
\]

The rank-0 representative is the zero tensor. The rank-1 minimal representative is the symmetric tensor with every entry equal to 0 except for the $(3,3,3,3)$th entry, which equals 1. The rank-2,3,4,5,6,7 minimal representatives are given respectively by 
\begin{align*}
&
\left[
\begin{array}{ccc|ccc|ccc}
0 & 0 & 0 & 0 & 0 & 0 & 0 & 0 & 0 \\
0 & 0 & 0 & 0 & 0 & 0 & 0 & 0 & 0 \\
0 & 0 & 0 & 0 & 0 & 0 & 0 & 0 & 0 \\ \hline
0 & 0 & 0 & 0 & 0 & 0 & 0 & 0 & 0 \\
0 & 0 & 0 & 0 & 0 & 1 & 0 & 1 & 1 \\
0 & 0 & 0 & 0 & 1 & 1 & 0 & 1 & 1 \\ \hline
0 & 0 & 0 & 0 & 0 & 0 & 0 & 0 & 0 \\
0 & 0 & 0 & 0 & 1 & 1 & 0 & 1 & 1 \\ 
0 & 0 & 0 & 0 & 1 & 1 & 0 & 1 & 1 \\ 
\end{array}
\right], 
\quad
\left[
\begin{array}{ccc|ccc|ccc}
0 & 0 & 0 & 0 & 0 & 0 & 0 & 0 & 0 \\
0 & 0 & 0 & 0 & 0 & 0 & 0 & 0 & 0 \\
0 & 0 & 0 & 0 & 0 & 0 & 0 & 0 & 0 \\ \hline
0 & 0 & 0 & 0 & 0 & 0 & 0 & 0 & 0 \\
0 & 0 & 0 & 0 & 0 & 1 & 0 & 1 & 1 \\
0 & 0 & 0 & 0 & 1 & 1 & 0 & 1 & 1 \\ \hline
0 & 0 & 0 & 0 & 0 & 0 & 0 & 0 & 0 \\
0 & 0 & 0 & 0 & 1 & 1 & 0 & 1 & 1 \\ 
0 & 0 & 0 & 0 & 1 & 1 & 0 & 1 & 0 \\ 
\end{array}
\right],
\\
&
\left[
\begin{array}{ccc|ccc|ccc}
0 & 0 & 0 & 0 & 0 & 1 & 0 & 1 & 0 \\
0 & 0 & 1 & 0 & 0 & 1 & 1 & 1 & 1 \\
0 & 1 & 0 & 1 & 1 & 1 & 0 & 1 & 0 \\ \hline
0 & 0 & 1 & 0 & 0 & 1 & 1 & 1 & 1 \\
0 & 0 & 1 & 0 & 0 & 1 & 1 & 1 & 1 \\
1 & 1 & 1 & 1 & 1 & 1 & 1 & 1 & 1 \\ \hline
0 & 1 & 0 & 1 & 1 & 1 & 0 & 1 & 0 \\
1 & 1 & 1 & 1 & 1 & 1 & 1 & 1 & 1 \\ 
0 & 1 & 0 & 1 & 1 & 1 & 0 & 1 & 0 \\
\end{array}
\right],
\quad
\left[
\begin{array}{ccc|ccc|ccc}
0 & 0 & 0 & 0 & 0 & 1 & 0 & 1 & 0 \\
0 & 0 & 1 & 0 & 0 & 1 & 1 & 1 & 1 \\
0 & 1 & 0 & 1 & 1 & 1 & 0 & 1 & 0 \\ \hline
0 & 0 & 1 & 0 & 0 & 1 & 1 & 1 & 1 \\
0 & 0 & 1 & 0 & 0 & 1 & 1 & 1 & 1 \\
1 & 1 & 1 & 1 & 1 & 1 & 1 & 1 & 1 \\ \hline
0 & 1 & 0 & 1 & 1 & 1 & 0 & 1 & 0 \\
1 & 1 & 1 & 1 & 1 & 1 & 1 & 1 & 1 \\ 
0 & 1 & 0 & 1 & 1 & 1 & 0 & 1 & 1 \\
\end{array}
\right],
\\
&
\left[
\begin{array}{ccc|ccc|ccc}
0 & 0 & 0 & 0 & 0 & 1 & 0 & 1 & 0 \\
0 & 0 & 1 & 0 & 0 & 1 & 1 & 1 & 1 \\
0 & 1 & 0 & 1 & 1 & 1 & 0 & 1 & 0 \\ \hline
0 & 0 & 1 & 0 & 0 & 1 & 1 & 1 & 1 \\
0 & 0 & 1 & 0 & 0 & 0 & 1 & 0 & 0 \\
1 & 1 & 1 & 1 & 0 & 0 & 1 & 0 & 0 \\ \hline
0 & 1 & 0 & 1 & 1 & 1 & 0 & 1 & 0 \\
1 & 1 & 1 & 1 & 0 & 0 & 1 & 0 & 0 \\ 
0 & 1 & 0 & 1 & 0 & 0 & 0 & 0 & 1 \\
\end{array}
\right],
\quad
\left[
\begin{array}{ccc|ccc|ccc}
0 & 0 & 0 & 0 & 0 & 1 & 0 & 1 & 0 \\
0 & 0 & 1 & 0 & 0 & 1 & 1 & 1 & 1 \\
0 & 1 & 0 & 1 & 1 & 1 & 0 & 1 & 0 \\ \hline
0 & 0 & 1 & 0 & 0 & 1 & 1 & 1 & 1 \\
0 & 0 & 1 & 0 & 0 & 0 & 1 & 0 & 0 \\
1 & 1 & 1 & 1 & 0 & 0 & 1 & 0 & 0 \\ \hline
0 & 1 & 0 & 1 & 1 & 1 & 0 & 1 & 0 \\
1 & 1 & 1 & 1 & 0 & 0 & 1 & 0 & 0 \\ 
0 & 1 & 0 & 1 & 0 & 0 & 0 & 0 & 0 \\
\end{array}
\right].
\end{align*}

We attempted to work over $\mathbb{F}_3$ but the computations could not be completed. We were able to determine that the maximum symmetric rank is at least 13 (since the program was terminated at this step while many symmetric tensors still had not been found). Thus, larger prime fields and larger tensor formats cannot be dealt with by our computer algorithm in a reasonable amount of time. 

\section{Conclusion}
We began by summarizing known results, and our own previous work about canonical forms of tensors (symmetric and non-symmetric). For third and fourth order tensors, the simplest cases for the classification of canonical forms is for $2 \times \dots \times 2$ tensor formats. The classification over $\mathbb{R}$ and $\mathbb{C}$ has appeared many times throughout the literature. In previous work, we extended the results by determining the canonical forms of $2 \times 2 \times 2$ tensors over the prime fields $\mathbb{F}_p$ for $p = 2, 3, 5$. We also examined the same problem for the larger tensor format $2 \times 2 \times 2 \times 2$ over $\mathbb{F}_p$ for $p = 2, 3$, where the choice of $p$ is restricted by the memory capabilities of our computer. 

Next, we summarized known results about $2 \times 2 \times 2$ symmetric and $2 \times 2 \times 2 \times 2$ symmetric tensors. The canonical forms of these tensor formats have been enumerated in the past over $\mathbb{R}$ and $\mathbb{C}$ by Weinberg \cite{Weinberg} and Gurevich \cite{Gurevich}. In recent work, we extended these results by determining the canonical forms over $\mathbb{F}_p$. For the third order case, we considered all primes $p \leq 17$. For the fourth order case, we considered $p \leq 5$. 

In this paper, we determined the canonical forms of second, third, and fourth order symmetric tensors of format $3 \times \dots \times 3$ over $\mathbb{F}_p$ for $p \in \{2, 3, 5\}$. We determined that the maximum symmetric rank over finite fields is at least 4 in the second order case, 7 in the third order case, and 13 in the fourth order case. Due to memory limitations we could not consider larger primes or tensor formats.

\section*{Acknowledgements}
The author would like to thank Edoardo Ballico for his insights and suggestions. 


\end{document}